\documentclass[11pt]{article}
\usepackage{graphicx}
\usepackage{algpseudocode}
\usepackage{algorithm}
\usepackage{amsmath}
\usepackage{amsfonts}
\usepackage[normalem]{ulem}
\usepackage{adjustbox}
\usepackage[shortlabels]{enumitem}
\usepackage{mathtools}
\usepackage{cite}
\usepackage{xcolor}
\usepackage{authblk}
\algnewcommand{\IIf}[1]{\State\algorithmicif\ #1\ \algorithmicthen}
\algnewcommand{\EndIIf}{\unskip\ \algorithmicend\ \algorithmicif}
\newcommand{\inner}[1]{\left\langle#1\right\rangle}

\def\R{\mathbb{R}}

\def\bone{{\mathbf 1}}

\def\ba{{\mathbf a}}
\def\bb{{\mathbf b}}

\def \bp {{\mathbf{p}}}
\def \bq {{\mathbf{q}}}

\def\bz{{\mathbf z}}

\def\bP{\mathbf P}
\def\diag{\mathbf{diag}}
\def\calP{\mathcal P}
\def\calA{\mathcal A}
\def\calB{\mathcal B}
\def\calE{\mathcal E}
\def\calL{\mathcal L}

\def\calM{\mathcal M}
\def\calN{\mathcal N}
\def\calS{\mathcal S}
\def\calT{\mathcal T}
\def\calB{\mathcal B}

\def\argmax{\mathop{\rm arg\,max}\limits}
\def\argmin{\mathop{\rm arg\,min}\limits}

\def\maxop{\mathop{\rm max}\limits} 

\def\Tbar{\mathbf{\overline{T}}}

\newcommand{\commentbyKG}[1]{}

\newcommand{\commentbyET}[1]{}
\newcommand{\alertbyET}[1]{{\color{blue} #1}}

\newtheorem{theorem}{Theorem}[section]
\newtheorem{lemma}[theorem]{Lemma}

\newtheorem{definition}[theorem]{Definition}

\numberwithin{equation}{section}
\date{}

\providecommand{\keywords}[1]
{
  \small	
  \textbf{\textit{Keywords---}} #1
}

\begin{document}
\title{A scalable solution for the extended multi-channel facility location problem}
%
%

\author[1]{Etika Agarwal}
\author[1]{Karthik S. Gurumoorthy}
\author[1]{Ankit Ajit Jain}
\author[1]{Shantala Manchenahally}
\affil[1]{Walmart Global Tech, Bangalore, India}
\affil[ ]{{\{etika.agarwal, karthik.gurumoorthy, ankit.ajit.jain, shantala.manchenahally\}@walmart.com}}



%
\maketitle 

\begin{abstract}
We study the extended version of the non-uniform, capacitated facility location problem with multiple fulfilment channels between the facilities and clients, each with their own channel capacities and service cost. Though the problem has been extensively studied in the literature, all the prior works assume a single channel of fulfilment, and the existing methods based on linear programming, primal-dual relationships, local search heuristics etc. do not scale for a large supply chain system involving millions of decision variables. Using the concepts of submodularity and optimal transport theory, we present a scalable algorithm for determining the set of facilities to be opened under a cardinality constraint. By introducing various schemes such as: (i) iterative facility selection using incremental gain, (ii) approximation of the linear program using novel multi-stage Sinkhorn iterations, (iii) creation of facilities one for each fulfilment channel etc., we develop a fast but a tight approximate solution, requiring $O\left(\frac{3+k}{2} m \ln\left(\frac{1}{\epsilon}\right)\right)$ instances of optimal transport problems to select $k$ facilities from $m$ options, each solvable in linear time. Our algorithm is implicitly endowed with all the theoretical guarantees enjoyed by submodular maximisation problems and the Sinkhorn distances. When compared against the state-of-the-art commercial MILP solvers, we obtain a 100-fold speedup in computation, while the difference in objective values lies within a narrow range of $\mathbf{3\%}$.\\\\
\end{abstract}
\section{Introduction}
Most supply chains in the online retail system can be represented as a bi-partite graph between sets of facilities and clients. Typically, each client has a demand $d_j$ to be serviced by one or many facilities. Each facility has a capacity $fcap_i$ representing the maximum amount of demand it can satisfy. Each facility is also associated with an one time fixed opening cost $F_{i}$ and can serve the demands only when it is open. The primary problem of interest is in finding the subset of facilities to be opened, and the assignment of fraction of demand serviced by the facility $i$ for the client $j$, so that the overall cost of opening facilities (fixed cost) and assigning demand (variable cost) is minimized. This is popularly known as the \emph{Facility Location Problem} (FLP) \cite{Frieze1974, Pal2001, Korupolu2000, Chudak2005} and is studied under different variants as discussed in Section~\ref{sec:relatedwork}. 

To the best of our knowledge, all the previous works on FLP assume a single channel of fulfilment ---a single edge of one type in the bi-partite graph---  between the facilities and clients which is used as a medium of service. The cost of servicing one unit of demand for the client $j$ by the facility $i$ using this one channel is denoted by the service cost $c_{ij}$. However, many real life scenarios exist where a facility can service a client via different channels, each with their own channel capacities and service cost. For instance, a source facility can send the required orders to the customers either through ground, 1-day air, or 2-day air shipments parcels etc. As opposed to a standard FLP, this introduces more complex constraints in the sense that the total outflow from a facility is not only limited by its capacity, but the capacities of individual shipment channels as well, leading to coupling constraints between the shipment channels.

In this paper, we study the extended version of non-uniform, capacitated facility location problem. We consider the scenario where: (i) the facility capacities are hard in that each facility can be opened at most once to serve a maximum demand of $fcap_i$, (ii) there are multiple fulfilment channels $ e \in E$ between the facilities and clients, where the cost incurred to serve one unit of demand by facility $i$ to client $j$ via the fulfilment channel $e$ is $c_{ije}$, (iii) each facility $i$ has a channel capacity $ccap_{ie}$ denoting the maximum demand fulfilled by it using $e$, across all clients, (iv) the maximum number of facilities that can be opened is limited to a user input $k$ – this is referred to as the cardinality constraint. This extended FLP can be formulated as an MILP with one Boolean variable for each facility denoting its open status, as described in Appendix~\ref{sec:MILP}. However, it is known that even in the most basic forms, a capacitated FLP is $\calN \calP$-hard \cite{Pal2001}, and the run time for the direct MILP solver grows exponentially in the number of facilities and hence is practically infeasible for a large-scale supply chain system.

\subsection{Contributions}
\label{sec:contribution}
We present a scalable solution framework for the extended FLP by formulating it as a set selection problem. FLP problems in this form are known to be submodular\cite{Frieze1974}. We utilize the stochastic distorted greedy algorithm proposed in \cite{Harshaw2019} for a special class of non-monotone submodular optimization problems, to arrive at a fast but good quality approximate solution for the FLP. The stochastic distorted greedy method incrementally selects the facility one at a time to maximize a distorted marginal gain. The calculation of best marginal gain requires repeatedly solving a demand allocation problem $\mathbb{P}_3$ in \eqref{eqn:demand_allocation} which is an instance of Linear Program (LP). In the context of greedy algorithm for our FLP, any algorithm used to solve $\mathbb{P}_3$ is referred to as the \textit{value oracle}.  The complexity of applying greedy method is in the development of a \textit{fast} value oracle. Though the LP can be solved in polynomial time \cite{Karmarkar1984, Khachiyan1979}, it does not scale for large-sized problems involving millions of decision variables, and more so when multiple instances of $\mathbb{P}_3$ need to be solved over many iterations. Our primary contribution is the development of a fast approximate algorithm for the value oracle to compute the solution of $\mathbb{P}_3$. We replace the direct LP solver by borrowing ideas from the domain of Optimal Transport (OT) theory \cite{peyre19a}. Though OT studies the problem of optimal allocation of resources between supply (facilities) and demand (client) locations, the presence of multiple fulfilment channels between them precludes the direct application of efficient OT algorithms proposed in~\cite{cuturi13a, altschuler2017, abid18a}, referred to as the Sinkhorn iterations. Notwithstanding this difficulty, we show how the LP involved in the value oracle can still be approximated by a multi-stage OT problem. We first merge the multiple shipment channels into a single abstract channel between every facility-client pair. The stage-1 OT problem is solved for this approximate network to optimally allocate the client demand to facilities. Subsequently, the shipment channels are split again and stage-2 OT is solved, one for each facility, to determine the distribution of allocated demand from stage-1 OT on to multiple shipment channels.  The second stage OT gives an approximate solution to the demand allocation problem $\mathbb{P}_3$ which is used in stochastic greedy algorithm. Although we incur an increase in the number of instances of OT problems compared to LP instances, each OT problem is computable in lightspeed using Sinkhorn iterations \cite{cuturi13a, altschuler2017, abid18a} and the overall computation time is much lower compared to executing a direct LP solver. Our work,
\begin{itemize}
    \item is the first to combine the strengths of both the submodularity and the optimal transport theory, towards solving FLP.
    \item introduces an extended version of FLP involving multiple fulfilment channels and presents a scalable algorithm for a large-scale supply chain system.
\end{itemize}

\subsection{Problem formulation}
\label{sec:problem}
Consider a supply chain network $\mathbf{T} = \left(\calM, \calN, E, \calE\right)$ with $n$ customers or clients, each with a demand $d_j$, $j \in \mathcal{N} = \{1, 2, \dots, n\}$, and $m$ facilities serving these demand nodes. Each facility $i$ has a one-time fixed cost to open $F_{i}$, $i \in \calM = \{1, 2, \dots, m\}$, associated with its operation. Facilities can supply to the demand destinations via multiple shipment (fulfilment) channels $e \in E$. Let $c_{i,j,e}$ be the cost to ship unit item from facility $i$ to client $j$ using the channel $e$. Each facility has a channel capacity $ccap_{i,e}$ denoting the maximum demand that facility $i$ can fulfill using the channel $e$, across all clients. The total outflow from a facility $i$ cannot exceed the maximum fulfilment capacity of $fcap_i$, $i \in \calM = \{1, 2, \dots, m\}$. Without the loss of generality, we assume $ccap_{i,e} \leq fcap_i$. It is not required that all the demand from a client be serviced by a single facility, and not all facilities can serve a client location using all the shipment channels pertaining to the practical limitations on transport connectivity. Let $\mathcal{E}$ be the set of all edge connections $( i, j, e )$ where a shipment channel $e $ connects facility $ i $ to client $ j $, and $\calN_i$ be the set of all clients that can be serviced by facility $i$. Denote by $E_i$, the shipment channels types available with the facility $i$.\commentbyET{\sout{For any $(i,j,e) \notin \mathcal{E}$, $c_{i,j,e} = \infty$} we need to make profit zero (infinite cost will not give zero profit).} The sub-network formed by the facility $i$, clients $j \in \calN_{i}$ and the shipment channels between them is represented by $\mathbf{T}_{i} = \left(\calN_{i}, E_{i}, \calE_{i}\right)$.

We aim to solve a facility location problem where the objective is to select an optimal set of facilities $\calS^* \subseteq \calM$, with a cardinality constraint $|\calS^*| \leq k $, to meet maximum demand with minimum costs incurred in shipment and opening new facilities, while adhering to the dynamics of the supply chain network. We consider a soft constraint on the demand, in the sense that it is allowed to fulfil only partial demand, but a penalty is incurred for any unfulfilled demand. Given the supply chain network described above, the optimal set of facilities $\mathcal{S}^\ast$ is obtained by solving the set selection problem $\mathbb{P}_1$ namely:
\begin{align} \label{eqn:set_selection1:objective}
 \mathbb{P}_1:    \mathcal{S}^\ast   = \argmin_{\calS\subseteq\calM, \quad |\calS|\leq k} &\Biggl(\sum_{i \in \calS} F_i + \min_{x_{i,j,e}} \biggl( \sum_{i \in \calS } \sum_{(i,j,e) \in \mathcal{E}_i} c_{i,j,e}  x_{i,j,e} d_{j} + \nonumber \\
 & C \Bigl( \sum_{j \in \calN} d_{j} - \sum_{i \in \calS } \sum_{(i,j,e) \in \calE_i} x_{i,j,e} d_{j} \Bigr)\biggr)\Biggr)
 \end{align}
 such that,
\begin{align}
\sum_{i \in \calS} \sum_{e \in E_i} x_{i,j,e} & \leq 1 & \quad \forall j \in \calN \label{eqn:set_selection1:demand} \\
\sum_{j \in \calN_i} \sum_{e \in E_i} x_{i,j,e} d_j & \leq fcap_i & \quad \forall i \in \calS \label{eqn:set_selection1:fcap}\\
\sum_{j \in \calN_i} x_{i,j,e} d_j &\leq ccap_{i,e} & \quad \forall i \in \calS, \quad \forall e \in E_i \label{eqn:set_selection1:ccap}\\
 x_{i,j,e} &\geq 0 & \quad \forall \left(i, j, e\right) \in \calE \label{eqn:set_selection1:positive_ship}\\
  x_{i,j,e} &\leq 0 & \quad \forall (i, j, e) \in \{\calS \times \calN \times E\}\backslash \cal{E}. \label{eqn:set_selection1:negative_ship}
  \end{align}    
The value $x_{i,j,e}$ is the proportion of demand $d_j$ at the $j^{th}$ client fulfilled by the $i^{th}$ facility using the shipment channel $e$. The constraints in \eqref{eqn:set_selection1:fcap} and \eqref{eqn:set_selection1:ccap} enforce the facility capacity and channel capacity constraints on the outflow from the selected facilities $\mathcal{S}$, and \eqref{eqn:set_selection1:positive_ship} and \eqref{eqn:set_selection1:negative_ship} ensure that the network connectivity constraints are followed. The cardinality constraint on the total number of selected facilities is imposed as $\left|\mathcal{S}\right|\le k$. The optimization objective of the FLP given in~\eqref{eqn:set_selection1:objective} has three components: (i) first term is the fixed cost to open the selected facilities $\mathcal{S}$, (ii) the middle term denotes the total cost to ship, and (iii) lastly a penalty on any unfulfilled demand. The parameter $C > \max\limits_{(i,j,e) \in \calE} c_{i,j,e}$ controls the trade-off between the fulfilled and unfulfilled demand. Some times an FLP solution may not be able to meet all the demand either due to insufficient capacity, or insufficient supply due to transportation restrictions (network constraints). Allowing a slack on the fulfilled demand in \eqref{eqn:set_selection1:demand} ensures the feasibility of the problem. For $\left(i,j,e\right)\in\mathcal{E}$, define $p_{i,j,e}=C-c_{i,j,e}\geq0$ as the \emph{profit} to ship unit item from facility $i$ to client $ j $ using the channel $e$, and $0$ otherwise. Using $p_{i,j,e}$ and dropping the constant $C\sum_{j\in\mathcal{N}} d_j$, the set selection problem $\mathbb{P}_1$ can be re-stated as a maximization problem over a set function,
\begin{align} \label{eqn:set_selection2}
    \mathbb{P}_2:\quad  \calS^\ast &= \argmax_{\calS\subseteq\calM, \quad |\calS|\leq k} \quad g(\calS) - h(\calS)
\end{align}
where $h\left(\mathcal{S}\right)=\sum_{i\in\mathcal{S}} F_i$ is the total cost of opening the facilities in $\calS$ and steadily increases as more facilities are included. The component $g\left(\mathcal{S}\right)$ is the optimal cost of demand allocation obtained by solving the sub-problem namely,
\begin{align} \label{eqn:demand_allocation}
    \mathbb{P}_3: g\left(\mathcal{S}\right)=\max_{x_{i,j,e}} {\sum_{i\in\mathcal{S}}\sum_{\left(i,j,e\right)\in\mathcal{E}_i}{p_{i,j,e}x_{i,j,e}d_j}}
\end{align}
subject to the constraints in \eqref{eqn:set_selection1:demand}-\eqref{eqn:set_selection1:negative_ship}. The objective here is to fulfill demand using the facilities $\mathcal{S}$ at the most profitable shipment value $p_{i,j,e}$. As any unfulfilled demand has profit $0$, it must be minimized for profit maximization in \eqref{eqn:demand_allocation}. It is easy to see that $g\left(.\right)\geq0$ and non-decreasing as for $\mathcal{S}_1\subseteq\mathcal{S}_2, g\left(\mathcal{S}_1\right)\le g\left(\mathcal{S}_2\right)$. This definition of FLP will allow us to utilize the concept of submodularity for a fast solution in the subsequent sections.

\section{Background}
\label{sec:background}

\subsection{Submodularity}
\label{sec:submodularity}
We briefly review the concept of submodular functions, which we later use to develop our algorithm. For a set function $f(\cdot)$, a subset $\calA \subseteq \calM $ and an element $i \in \calM$, define the incremental gain as $f_{\calA}(i) = f\left(\calA \cup \{i\}\right) - f(\calA)$.
\begin{definition}[Modularity, Submodularity and Monotonicity]
\label{def:sub}
Consider any two sets $\calA \subseteq \calB \subseteq \calM$. A set function $f(\cdot)$ is submodular iff for any $i \notin \calB$, $  f_{\calA}(i) \geq f_{\calB}(i)$. It is called monotone when $f(\calA)\leq f(\calB)$ and modular when $f_{\calA}(i)$ is independent of $\calA$.
\end{definition}
Submodularity implies diminishing returns where the incremental gain in adding a new element $i$ to a set $\calA$ is at least as high as adding to its superset $\calB$~\cite{fujishige05, lo83}. 
\commentbyKG{
Another characterization of submodularity is via the submodularity ratio \cite{weaksub,weaksubInit} defined as follows.
\begin{definition}[Submodularity Ratio]
Given two disjoint sets $\calL$ and $\calS$, and a set function $f(\cdot)$, the submodularity ratio of $f(\cdot)$ for the ordered pair ($\calL,\calS$) is given by:
\begin{equation}\label{eqn:subratio}
\alpha_{\calL,\calS} \coloneqq \frac{\sum\limits_{i \in \calS} f_{\calL}(i)}{f\left(\calL \cup \calS\right)-f(\calL)}.
\end{equation}
\end{definition}
\alertbyET{Do we need Definition 2.2 and this paragraph?}Submodularity ratio captures the increment in $f(\cdot)$ by adding the entire subset $\calS$ to $\calL$, compared to summed gain of adding its elements individually to $\calL$. It is known that $f(\cdot)$ is submodular if and only if $\alpha_{\calL,\calS} \geq 1, \forall \calL,\calS$. In the case where $\alpha_{\calL,\calS} <1$ but bounded away from $0$, $f(\cdot)$ is called \textit{weakly submodular}~\cite{weaksubInit}. 
}
Submodular functions enjoy provable, tight performance bounds when the set elements are selected incrementally and greedily~\cite{Nemhauser78,weaksub,proto}.

\subsection{Optimal transport (OT)}\label{sec:OT}
Let $\calM = \{1, 2, \dots, m\}$ and $\calN=\{1, 2, \dots, n\}$ be a set of supply and demand locations, associated with a capacity $\bp_i$ and demand $\bq_j$ at locations $i \in \calM$ and $j\in \calN$ respectively.\commentbyKG{ The marginals $\bp$ and $\bq$ lie on corresponding simplices $\Delta_{m}$ and $\Delta_n$, where $\Delta_{k}\coloneqq\{\bz\in\R_{+}^k|\sum_i \bz_i=1\}$.} The OT problem \cite{kantorovich42a} aims at finding a most profitable transport plan $\gamma$ as a solution to\footnote{When formulated as a minimisation problem, the objective is to find a transport plan that yields the minimal transporting effort.} $\maxop_{\gamma\in\Gamma(\bp,\bq)} \inner{\bP,\gamma} $,
where $\Gamma(\bp,\bq)\coloneqq\{\gamma\in\R_{+}^{m\times n}|\gamma\bone \leq \bp;\gamma^\top\bone \leq \bq\}$\commentbyET{is the space of joint measure between the supply and the demand marginals, $\bp$ and $\bq$ respectively}. Here, $\bP\in\R_{+}^{m\times n}$ with entries $\bP_{ij}$ representing the \textit{profit} of transporting a unit mass from supply $i \in\calM$ to demand $j\in\calN$. Recently, \cite{cuturi13a} proposed an extremely fast solution for learning entropy regularized transport plan $\gamma$ using the Sinkhorn algorithm~\cite{knight08a}. For a recent survey on OT, please refer to \cite{peyre19a}. 

\section{A scalable approximate solution}
\label{sec:scalablealgo}
Facility location problems are known to be submodular \cite{Frieze1974, Chudak2005}. In the context of FLP, the diminishing returns property of sub-modularity implies that it is more beneficial to open a new facility when there are fewer facilities already open. Maximisation of submodular set functions is a well studied problem. As they are $\calN \calP$-hard in the general form, a variety of algorithms have been proposed to find approximate solutions to submodular optimization problems. One of the most popular category of algorithms are the variants of incremental selection of a set using greedy approaches~\cite{Nemhauser78,Buchbinder2014, Buchbinder2017, Kuhnle2019, Sakaue2020}.
 These algorithms  provide an approximation guarantee to the optimal solution if and only if the submodular objective is monotone and non-negative. However, in our case, while both $g\left(\mathcal{S}\right)$ and $h\left(\mathcal{S}\right)$ are monotone, the overall optimization objective in $\mathbb{P}_2$ is non-monotone and could also be negative. Therefore, a vanilla greedy algorithm does not give any guarantees on the quality of the solution. We note that the set function $g\left(\mathcal{S}\right)$, obtained as the output of the optimal demand allocation in $\mathbb{P}_3$ is non-negative monotone submodular function, and $h\left(\mathcal{S}\right)$ is monotone modular. This allows us to use a more recent \textit{distorted greedy} algorithm proposed in \cite{Harshaw2019}. 
 
 A faster variation of the same algorithm called the \textit{stochastic distorted greedy} is outlined in Algorithm~\ref{algo:distorted_greedy}, and progresses as follows. Let $\mathcal{S}$, $|\calS| = \ell - 1$, represent the set of currently chosen facilities. The property of diminishing returns enables us to select facilities in an incremental fashion, where starting from the empty set $\mathcal{S}=\emptyset$, in every iteration we identify the next best facility $u_\ell$ from a random subset $\mathcal{M}^r$ of available facility options $\mathcal{M}-\mathcal{S}$. The selection is made to maximize a modified incremental gain $\left(1-\frac{1}{k}\right)^{k-\ell} g_{\calS}(u_\ell) - h_{\calS}(u_\ell)$, where $\left(1-\frac{1}{k}\right)^{k-\ell}$ is a continuously changing distortion factor applied to $g_\calS(u_\ell) = g\left(\calS \cup \{u_\ell\}\right) - g(\calS)$. The set $\mathcal{S}$ is grown incrementally by including the chosen facility $u_\ell$ and the entire process is repeated for $k$ iterations. As more and more facilities are added, the distortion factor puts increasingly more emphasis on the reduction of shipment cost by appropriately scaling $g(\calS)$. This way of incrementally selecting the facility set eliminates the need of integrality constraints and is a major factor for reducing the solution complexity. The recommended size of the random set $\calM^r$ is $r  = \left \lceil \frac{m}{k} \log \left( \frac{1}{\epsilon} \right) \right \rceil$, where $\epsilon$ is a hyper-parameter which directly affects the approximation guarantees. The total number of invocations to set function $g(\cdot)$ over all the $k$ iterations is $O\left(m \ln\left(\frac{1}{\epsilon}\right)\right)$, and is independent of $k$. Let $\calS$ be the approximate solution from stochastic distorted greedy algorithm, and $S^{\ast}$ the optimal solution. The stochastic distorted greedy algorithm gives guarantees on the optimality gap in an expected sense namely, $\mathbb{E}\left[g(\calS) - h(\calS)\right] \geq \left(1-\frac{1}{e}-\epsilon\right)g\left(S^{\ast}\right) - h\left(S^{\ast}\right).$ If all facilities have the same cost to open $(F_{i})$, the approximation guarantee match that in \cite{mirzasoleiman15a} for monotone sub-modular functions. 
\begin{algorithm}[t]
\caption{Stochastic distorted greedy algorithm}\label{algo:distorted_greedy}
{
\begin{algorithmic}[1]
\State {Set $r = \left\lceil \frac{m}{k} \log\left(\frac{1}{\epsilon} \right) \right\rceil$, $\calS = \emptyset$}
\For{$\ell = 1$ to $k$}
    \State $\calM^r \gets$ select $r$ elements randomly from $\calM-\calS$
    \For{$u \in \calM^r$} 
        \State  Calculate $g\left(\calS \cup \{u\}\right)$, $g_{\calS}(u) \gets g\left(\calS\cup \{u\}\right) - g\left(\calS\right)$, $h_{\calS}(u) \gets f_u$
    \EndFor
    \State $u_\ell \gets \argmax_{u \in \calM^r} \left\{ \left(1-\frac{1}{k}\right)^{k-\ell} g_{\calS}(u) - h_{\calS}(u)\right \}$
    \State $IG \gets \left(1-\frac{1}{k}\right)^{k-\ell} g_{\calS}(u_\ell) - h_{\calS}(u_\ell)$
    \IIf{$IG > 0$} {$\calS \gets \calS \cup \{u_\ell\}$} 
\EndFor
\end{algorithmic}}
\end{algorithm}

\subsection{Fast value oracle}
The greedy algorithm and its variants for the optimal set selection problem measure the computational complexity of the solution in terms of the number of calls to a \textit{value oracle} \cite{Harshaw2019, Buchbinder2017, Kuhnle2019, Sakaue2020}. Given a set of facilities $\calS$, a \textit{value oracle} is a system or an algorithm to compute the value of the submodular function. As the greedy algorithm and its variants require repeated evaluation of the incremental gain, \textit{value oracle} is used multiple times in the set selection process. For the solution of FLP problem in $\mathbb{P}_2$ using the stochastic distorted greedy method in Algorithm~\ref{algo:distorted_greedy}, the \textit{value oracle} must solve the optimal demand allocation problem $\mathbb{P}_3$ for $O\left(m \ln\left(\frac{1}{\epsilon}\right)\right)$ different instance of LP, and hence can take significant amount of time for a large scale supply chain network involving millions of decision variables. In this section, we propose to solve the demand allocation problem using Sinkhorn iterations~\cite{cuturi13a}, popular in optimal transportation theory, for a fast, but good quality solution. As our problem is not a standard OT, we start by introducing methods to reduce the demand allocation problem $\mathbb{P}_3$ into a novel multi-stage optimal transport problem.

\subsubsection{Channel decoupling under favourable conditions}\label{sec:channel_decoupling}
The presence of multiple fulfilment channels in the problem $\mathbb{P}_3$ introduces complexity by coupling the channels through the shared facility capacity $fcap_i$. If $fcap_i$ was large enough for all the facilities, it would no longer be a limiting constraint, and every facility and fulfilment channel pair could be treated independently. The limit over which the facility capacity in \eqref{eqn:set_selection1:fcap} plays no role to the demand allocation problem $\mathbb{P}_3$ is given in Lemma \ref{lemma:decouple1}. Proofs are available in Appendix~\ref{sec:prooofs}.

\begin{lemma} \label{lemma:decouple1}
The facility capacity constraint is superfluous and trivially satisfied for all facilities $i$ satisfying $\sum_{e\in E_i} ccap_{i,e} \leq fcap_i$. 
\end{lemma}
It may further be possible to decouple the fulfilment channels of facilities with not enough facility capacity using the result below. 
\begin{lemma}\label{lemma:decouple2}
Consider a facility $i \in \calS$ such that $\sum_{e \in E_i} ccap_{i,e} > fcap_i$. If $\exists e^{\prime} \in E_i$ satisfying
$
    e^{\prime} = \argmin\limits_{e \in E_i} p_{i,j,e}, \mbox{\hspace{5pt}} \forall j \in \calN_i, \label{eqn:costly_edge}
$,
the channel capacity of $e^{\prime}$ at the facility $i$ can be reduced to: $ccap_{i,e^{\prime}} = \max\left(0,fcap_i - \left(\sum_{e \in E_i-\{e^{\prime}\}} ccap_{i,e}\right)\right)$.
\end{lemma}
Lemma \ref{lemma:decouple2} states that we can reduce the capacity requirement of the least profitable fulfilment channel $e^{\prime}$ by modifying its \textit{ccap}. \commentbyET{Note that if a client $j$ is not serviced by some channel $e \in E_i$ at facility $i$,  the transportation profit $p_{i,j,e^{\prime}}$ is set to $0$ and the least profitable fulfilment channel $e^{\prime}$ is chosen accordingly.} However, this reduction is possible only when all the clients serviced by $e^{\prime}$ are also serviced by all the other channels. Note that for clients not serviced by $e^{\prime}$, the condition in Lemma \ref{lemma:decouple2} is trivially satisfied as $p_{i,j,e^{\prime}} = 0$. If the facility capacity is still not enough, that is, $\sum_{e \in E_i} ccap_{i,e} > fcap_i$ even after the least profitable channel capacity reduction, the process can be repeated by applying Lemma~\ref{lemma:decouple2} on the next least profitable channel-type with non-zero channel capacity satisfying \eqref{eqn:costly_edge}. This way, a repeated application of Lemma \ref{lemma:decouple2}, and lastly an application of Lemma \ref{lemma:decouple1} can be used to decouple the fulfilment channels. If the fulfilment channels at facility $i$ can be completely decoupled this way, then specific to the problem $\mathbb{P}_3$, each fulfilment channel associated with a facility $i$ can be treated as a separate facility indexed by (hyphenated letters) $i{\text -}e$, of capacity $fcap_{i{\text -}e}= ccap_{i,e}$, with profits on the \emph{single} edge $e$ to client $j$ given by $p_{i{\text -}e, j}=p_{i,j,e}$. We explain the channel decoupling process with an example in Appendix~\ref{sec:channeldecouplingexample}. Henceforth in the paper we assume that the channel decoupling is already applied for problem $\mathbb{P}_3$ wherever possible, and for the sake of brevity continue to use the same original notations as before to represent this partially decoupled network model as well.

\subsubsection{Reduced channel model}\label{sec:reduced_model}
\commentbyET{Part of the results from Section~\ref{sec:channel_decoupling} can be used to decouple the shipment channels and make them independent of the outflow on other channels. Specifically, the facilities in the set $\calM_t$ satisfying Lemma~\ref{lemma:decouple1} ---after repeated application of Lemma~\ref{lemma:decouple2}--- can be expanded to $\calM_t^E = \{i{\text -}e: i \in \calM_t, e \in E_i\}$ such that, every facility in $\calM_t^E$ is connected to its respective clients using only a single fulfilment channel with its corresponding capacity $fcap_{i{\text -}e}$ and profits $p_{i{\text -}e, j}$. }

Results from the last section can be used to decouple the shipment channels and make them independent of the outflow on other channels. However, such straightforward reductions are not always possible when the conditions in Lemmas~\ref{lemma:decouple1} and \ref{lemma:decouple2} are not met. In this section we propose approximations for the remaining set of facilities with the non-trivial facility capacity constraint. Consider any facility $i \in \calM$. We introduce a new abstract shipment channel $\alpha$ and replace multiple shipment channels $E_i$ at the facility $i$ with $\alpha$\commentbyKG{ (see Fig. \ref{fig:edge_roll_up} for reference)} to get an approximate network model $\Tbar = \left(\calM, \calN, \overline{\calE}\right)$ comprised of equivalent facilities with only single channel of shipment $E_i = \alpha$, $\forall i \in \calM$, and edges $\overline{\calE} = \left\{(i,j,\alpha): i \in \calM \text{ and } j \in \calN_i\right\}$. Since there is only one channel $\alpha$, it is superfluous to impose its channel capacity as $ccap_{i,\alpha}=fcap_i$. The profit to ship one unit of item from the facility $i$ to a client $j$ through $\alpha$, is computed as a weighted sum of the profits on the original shipment channels $E_i$, namely
$p_{i,j,\alpha} := \sum_{e \in E_i} w_{i,e}p_{i,j,e}, \hspace{5pt} \forall j \in \calN_i$, where $w_{i,e} = \frac{ccap_{i,e}}{\sum_{e\in E_i} ccap_{i,e}}$ $\forall i\in\calM, e\in E_i$ represents the contribution of a channel to the total available channel capacity and is independent of the clients $j \in \calN_i$. This definition ensures that the profit $p_{i,j,\alpha}$ on the abstract channel, is balanced between both the original channel profits $p_{i,j,e}$ and its capacity contribution proportion $w_{i,e}$. As the network $\Tbar$ only comprises of the facility capacities, the problem reduces into an instance of standard OT discussed in Section~\ref{sec:OT}.

\subsubsection{Multi-stage Sinkhorn algorithm}
As discussed before, a quick solution to the FLP using the stochastic distorted greedy method would require \textit{fast value oracle} to solve the demand allocation problem $\mathbb{P}_3$, especially when it is repeatedly invoked over several iterations. The optimal transportation literature~\cite{peyre19a} provides a linear time algorithm to compute approximate but good quality solution using Sinkhorn iterations~\cite{cuturi13a,abid18a} for LP instances arising in this field. We identify that the same algorithm can be applied to demand allocation problem for a supply chain network with single shipment channel. By adding an entropy maximization term $(-\mu \sum\limits_{i,j} \gamma_{i,j} \log\left(\gamma_{i,j}\right))$ to the objective, the transport plan can be computed through an iterative algorithm outlined in Algorithm~\ref{algo:sinkhorn} involving simple matrix operations. While the solution to the OT lies on one of the vertex of the polytope that forms the feasible region, the addition of entropy term shifts it to the interior of the polytope. In the context of supply chain, this implies that a solution using Sinkhorn iterations is dense in terms of the used facility-client connections. The trade-off between sparsity and speed is controlled using a parameter $\mu>0$. Building on this framework, we propose a novel multi-stage Sinkhorn algorithm as \emph{the} fast value oracle for the multi-channel FLP.

Let $\Tbar_\calS$ be the sub-network from $\Tbar$ formed of the facilities in $\calS$, their clients, and the single shipment channel $\alpha$. The demand allocation problem in $\mathbb{P}_3$ is divided into two stages, namely demand allocation in (i) $\Tbar_\calS$ and (ii) $\mathbf{T_i}$, $\forall i \in \calS$, each modeled as an instance of OT. The demand allocation problem for $\Tbar$ with the selected set of facilities $\calS$ can be posed as the OT below,
\begin{equation}\label{eqn:OT_Tbar}
\mathbb{P}_4:   \quad \maxop_{\overline{\gamma}\in\Gamma(\overline{\bp},\overline{\bq})} \inner{\overline{\bP}_{\calS},\overline{\gamma}},
\end{equation} 
where $\overline{\bp}$ is the vector of facility capacities in $\calS$, and $\overline{\bq}$ is the vector of client demands $d_j$. $\overline{\bP}_{\calS}$ is the profit matrix with entries $\overline{\bP}_{\calS_{i,j, \alpha}} = p_{i,j,\alpha}$ on the edges in $\overline{\calE}$ corresponding to the facilities in $\calS$\commentbyET{, and $0$ otherwise}. We obtain the solution for $\mathbb{P}_4$ using the Sinkhorn iterations~\cite{cuturi13a} outlined in Algorithm~\ref{algo:sinkhorn}. It is adapted for a maximisation OT problem and differs slightly from the original algorithm proposed in~\cite{cuturi13a}. The solution of this \textit{stage-1 OT} problem is the desired proportion $\overline{x}_{i,j} = \frac{\overline{\gamma}_{i,j}}{d_j}$ of the demand $d_j$, to be fulfilled by the facilities $i \in {\calS}$.\commentbyET{This is in fact the recommended proportion of demand $d_j$ to be fulfilled by facility $i$.} The distribution of this demand among different fulfilment channels $x_{i,j,e}$ is not known yet.

\begin{algorithm}[!t]
\caption{Sinkhorn iterations $(\bp, \bq, \bP)$}\label{algo:sinkhorn}
{
\begin{algorithmic}[1]
\State If $\bone^T \bp < \bone^T \bq$, create a pseudo-supply node $\tilde{i}$ with capacity $\bp_{\tilde{i}} = \bone^T \bq -\bone^T \bp$,
\Statex  and append it to $\bp$. Set profit $P_{\tilde{i},j} = 0, \forall j$.
\State If $\bone^T \bq < \bone^T \bp$, create a pseudo-demand node $\tilde{j}$ with demand $\bq_{\tilde{j}} = \bone^T \bp - \bone^T \bq$,
\Statex  and append it to $\bq$. Set profit $P_{i,\tilde{j}} = 0, \forall i$.
\State Let $\kappa = \bone^T \bp = \bone^T \bq$. Define $\tilde{\bp}=\bp/\kappa$, $\tilde{\bq}=\bq/\kappa$.
\State Construct $\Omega$ with entries $ \Omega_{ij} = \exp\left(\frac{\bP_{ij}}{\mu}\right)$.
\State \textbf{Initialize:} $\ba=\tilde{\bp}$, $\bb = \tilde{\bq}$.
\While{not convergence}
    \State $\ba = \Big(\diag\big(\Omega \bb \big)\Big)^{-1} \tilde{\bp}$
    \State $\bb = \Big(\diag\big(\Omega^T \ba \big)\Big)^{-1} \tilde{\bq}$
\EndWhile
\State \textbf{Return:} $\gamma = \kappa*\left(\diag(\ba)\right) \mathbf{\Omega} \left(\diag(\bb)\right)$
\end{algorithmic}}
\end{algorithm}

\begin{algorithm}[!t]\label{algo:multi_sinkhorn}
\caption{Multi-stage Sinkhorn $(\calS)$}\label{algo:dual_sinkhorn}
{
\begin{algorithmic}[1]
\commentbyET{\State $\overline{\calM}_{\calS} = \{i{\text -}\alpha: i \in \calS\}$  \Comment{Abstract facilities corresponding to $\calS$}}
\State $ \overline{\bq}_j = d_j$, $\forall j \in \calN$\Comment{Demand at the client locations}
\State $\overline{\bp}_{i} = fcap_{i}$, $\forall i \in \calS$  \commentbyET{\Comment{Abstract facility capacities}}
\State $\overline{\bP}_{\calS_{i, j , \alpha}} = p_{i,j,\alpha}$, $\forall i \in \calS$, $j \in \calN$ \Comment{Profits for abstract single channel $\alpha$}
\State Solve $\mathbb{P}_4$ to get $\overline{\gamma}$: \textbf{run Algorithm \ref{algo:sinkhorn}} with $\left(\overline{\bp}, \overline{\bq}, \overline{\bP}_\calS\right)$
\For{$i \in \calS$}
    \commentbyET{\If{$i \in \calM_\alpha$}}
        \State $\bq^i_j = \overline{\gamma}_{i,j}$, $\forall j \in \calN_i$ \Comment{Use output of Step 4 here}
        \State $\bp^i_e = ccap_{i,e}$, $\forall e \in E_i$
        \State $\bP^i_{e,j} = p_{i,j,e}$, $\forall e \in E_i$, $j \in \calN_i$ 
        \State Solve $\mathbb{P}_5^i$ to get $\gamma^i$: \textbf{run Algorithm \ref{algo:sinkhorn}} with $\left(\bp^i, \bq^i, \bP^i\right)$
        \State $x_{i,j,e} = \frac{\gamma^i_{e,j}}{\overline{\gamma}_{i,j}}$
    \commentbyET{\Else
         \State $x_{i,j,e} = \overline{\gamma}_{i{\text-}e,j}/d_j$, $\forall e \in E_i$, $j \in \calN_i$
    \EndIf}
\EndFor
\end{algorithmic}}
\end{algorithm}
  A second round of Sinkhorn iterations are applied to each sub-network $\mathbf{T_i}$, $\forall i \in \calS$ to optimally distribute the allocated demand portion $x_{i,j,\alpha}$ among multiple fulfilment channels which were abstracted using the channel $\alpha$. The \textit{stage-2 optimal transportation} problem is
  \begin{equation}\label{eqn:OT_Ti}
\mathbb{P}_5^i: \maxop_{{\gamma^i}\in\Gamma(\bp^i,\bq^i)} \inner{\bP^i,\gamma^i},
\end{equation} 
where $\bp^i \in \R_{+}^{E_i}$ is the vector of channel capacities with $\bp^i_e = ccap_{i,e}$, and $\bq^i \in \R_{+}^{\calN_i}$  is the stage-1 OT solution vector at the $i$-th facility with $\bq^i_j = \overline{\gamma}_{i,j}$. The matrix $\bP^i \in\R_{+}^{E_i \times \calN_i}$ is the profit matrix with values $\bP^i_{e,j} = p_{i,j,e}$. From the solution $\gamma^i$, the channel level distribution can be computed as: $ x_{i,j,e} = \frac{\gamma^i_{e,j}}{\overline{\gamma}_{i,j}}$.

\subsubsection{Computational complexity}
The multi-stage Sinkhorn method is outlined in Algorithm~\ref{algo:dual_sinkhorn}. When used as the value oracle in the stochastic distorted greedy technique described in Section~\ref{sec:scalablealgo}, it requires utmost $\sum_{i=1}^k (i+1)r = \frac{3+k}{2} m \ln\left(\frac{1}{\epsilon}\right)$ calls to the Sinkhorn method in~Algorithm~\ref{algo:sinkhorn}, where $r=\left\lceil \frac{m}{k} \log\left(\frac{1}{\epsilon} \right) \right\rceil$ is the cardinality of the randomly selected set $\mathcal{M}^r \subseteq \mathcal{M}-\mathcal{S}$ in each iteration of stochastic distorted greedy algorithm. In particular, only $m \ln\left(\frac{1}{\epsilon}\right)$ invocations to Algorithm~\ref{algo:sinkhorn} are required in the stage-1 for the problem $\mathbb{P}_4$ in~\eqref{eqn:OT_Tbar}. A majority $\left(\frac{1+k}{2} m \ln\left(\frac{1}{\epsilon}\right) \right)$ calls are required only in the stage-2 OT for the problem $\mathbb{P}_5$ which are: (i) executed on the much smaller single facility sub-networks $\mathbf{T_i}$ and, (ii) can be computed in parallel across all the facilities. This results in an extremely fast value oracle to determine an approximate solution to the problem $\mathbb{P}_3$ in~\eqref{eqn:demand_allocation}.

\section{Related work}
\label{sec:relatedwork}
Multiple variants of FLP are studied in the literature. In the uncapacitated FLP (UFLP), the facility capacity is assumed to be $\infty$ \cite{Cornuejols1990, Shmoys1997}. In the capacitated setting (CFLP), the capacities are finite and upper bound the quantum of demand that a facility can serve. Some papers investigate the specialised case of uniform facility capacities where $fcap_i = fcap$, $\forall i$ \cite{Korupolu2000, Chudak2005}, while others consider the general scenario of non-uniform facility capacities \cite{Pal2001, Barahona1998}. FLPs are also classified based on whether the capacities are \emph{hard} or \emph{soft}. In case of hard capacities, each facility $i$ can be opened at most once and can serve a demand of at most $fcap_i$ \cite{Pal2001}, whereas the soft capacitated facilities can service demands marginally exceeding their $fcap_i$ but with an added penalty~\cite{Arya2001, Chudak1999, Jain1999, Charikar1999}. Another FLP variant allows facilities to be opened $k > 1$ times by replicating itself \cite{Shmoys1997, Chudak1999, Mahdian2003, Chudak2005}, and service demand up to $k \ast fcap_i$ at an opening cost of $k F_{i}$.

\commentbyET{Solving for the Boolean variables $y_i$ denoting the open status of facilities makes FLP $\mathcal{NP}$-hard.} As solving the FLP directly as a MILP is $\mathcal{NP}$-hard~\cite{Lund1994}, plethora of approximation algorithms based on techniques such as LP rounding, primal-dual algorithms, and local search heuristics have been developed to yield constant factor approximations. LP techniques are used in~\cite{Shmoys1997} to give the first $3.16$ approximation algorithm for the UFLP, and subsequently are extended in~\cite{Jain1999} to handle soft capacitated FLPs. While the work in~\cite{Chudak1999} develop an approximation algorithm for CFLP with uniform, soft capacities, \cite{Arya2001} present a local search algorithm for the nonuniform, soft capacity variant. Local search heuristics have been the most popular approach for hard optimization problems. For the CFLP, \cite{Korupolu2000} gave the first iterative approximation algorithm of value no more than $8+\epsilon$ times optimum, for any $\epsilon > 0$ by using the simple local search heuristic originally proposed in~\cite{Kuehn1963}. By using the property of submodularity, \cite{Chudak2005} improved on the analysis of the same heuristic and proved a stronger $6(1+\epsilon)$ approximation. 

For the nonuniform CLP variant with hard capacities considered in this paper, \cite{Pal2001} developed a approximation algorithm with a guarantee of $9+\epsilon$ under the assumption of a single channel of fulfilment between facilities and client, unlike the extended multi-channel setup which is the focus of our work. Though the technique in~\cite{Pal2001} has a polynomial computational complexity, it is not designed for a large-scale network as it involves many operations in each iteration such as: (i) add a new facility to the current solution $\calS$: add$(i)$, (ii) open one facility $i$ and close a subset of facilities $\calT \subseteq \calS-\{i\}$: open$\left(i,\calT\right)$, (iii) close one facility $i\in \calS$ and open a subset of facilities $\calT \subseteq \calM-\{i\}$: close$\left(i,\calT\right)$. The operations open$\left(i,\calT\right)$ and close$\left(i,\calT\right)$ are computationally demanding as they require choosing a subset of facilities $\calT$. Through an involved sequence of analysis, the authors prove that their algorithm terminates after a polynomial number of these operations to yield a $9+\epsilon$ approximation, and for each choice of $i$ identifies a set $\calT$ in polynomial time to execute open$\left(i,\calT\right)$ and close$\left(i,\calT\right)$. Each of these operations require solving an instance of a LP, similar to the problem $\mathbb{P}_3$ in~\eqref{eqn:demand_allocation}, referred here as the value oracle.

For a large-scale supply network running a large number of LP instances, though polynomial in number, will not scale. As explained in section~\ref{sec:contribution}, our contribution is not only to study the extended nonuniform CFLP but to also design an efficient implementation of the value oracle, by solving these LPs via the proposed lightspeed \emph{multi-stage} Sinkhorn iterations using concepts from OT~\cite{cuturi13a, altschuler2017, abid18a, dvurechensky18a, xie2022accelerated}. Leveraging the recent advancements in submodular optimisation~\cite{Buchbinder2014, Buchbinder2017, Kuhnle2019, Sakaue2020, Harshaw2019}, we develop a scalable algorithm by incrementally selecting the facilities using the stochastic distorted greedy method~\cite{Harshaw2019}, involving only the add$(i)$ operation. Our approach does not require opening or closing subsets of facilities in each iteration as done in~\cite{Pal2001}. The approximation guarantee of $1-\frac{1}{e} -\epsilon$ for the stochastic distorted greedy algorithm ensures that our solution will be closer to the global optimum.

\section{Experiments} \label{sec:experiments}

We first present results on a mid-size e-commerce network with $n=2000$ client locations, $m=150$ tentative facilities, and $|E| = 3$ different shipment channels\commentbyKG{\footnote{Our non-disclosure agreement precludes us from making the code publicly available.}}. The network is fairly dense with a total of over 800,000 flow paths between the facility-client pairs, each with their own cost to ship per unit quantity. There is sufficient variation in the cost of opening the facilities with the ratio of std. dev. to its mean being $30\%$. The corresponding ratios for $fcap$ and $ccap$ are $28\%$ and $24\%$ respectively. The cost penalty for unfulfilled client demands is set to $C = 5\times\max(c_{i,j,e})$. We solve the FLP problem for this network for different values of cardinality constraint $k$ using 4 different approaches: 
\begin{enumerate}[(i)]
    \item \emph{Method 1}: A direct MILP solver,
    \item The stochastic distorted greedy method in Algorithm~\ref{algo:distorted_greedy} with $\epsilon = 0.01$ and following three different value oracles,
    \begin{enumerate}
        \item \emph{Method 2}: a standard commercial LP solver,
        \item \emph{Method 3}: the proposed multi-stage Sinkhorn algorithm,
        \item \emph{Method 4}: single-stage Sinkhorn algorithm.
    \end{enumerate}
\end{enumerate}
\commentbyET{As we directly compare our solution with the optimal solution from MILP to gauge the accuracy and quality of our approach, we find it redundant to baseline against other algorithms such as~\cite{Pal2001} which are computationally intensive.}

The \emph{single-stage Sinkhorn} oracle in Method 4 outputs the objective value of stage-1 OT problem in~\eqref{eqn:OT_Tbar}, invoking the Algorithm~\ref{algo:sinkhorn} only for stage-1. This is equivalent to selecting all the facilities based on the reduced network $\Tbar = \left(\calM, \calN, \overline{\calE}\right)$ constructed in Section~\ref{sec:reduced_model}. This oracle will incur the least computation time requiring only $m \ln\left(\frac{1}{\epsilon}\right)$ (independent of $k$) invocations to Algorithm~\ref{algo:sinkhorn}. Though it cannot be used to obtain a solution for the problem $\mathbb{P}_3$ in~\eqref{eqn:demand_allocation} as the demand distribution among different fulfilment channels is not computed, we used it as a test bed to also assess the quality of our reduced network $\Tbar$.
\begin{figure}[t] 
  \centering
  \includegraphics[width=0.9\linewidth,trim=0cm 0cm 0cm 0cm]{results_small_net.png}
  \caption{Results for different $k$ values. (a) Time taken to select facilities, (b) percentage overlap between selected facilities, (c) relative difference between the objective values in \eqref{eqn:set_selection1:objective}, and (d) relative difference in objective \eqref{eqn:demand_allocation} calculated with different value oracles.}
  \label{fig:compare}
\end{figure}
 After the facilities are selected, the demand allocation problem $\mathbb{P}_3$ is finally solved once using a standard LP solver with the final set of recommended facilities from these algorithms, and the objective value $J$ is calculated as per equation~\eqref{eqn:set_selection1:objective} for a fair comparison between approaches. Let $\mathcal{S}_{MILP}$, $\mathcal{S}_{LP}$, $\mathcal{S}_{Sink_2}$ and $\mathcal{S}_{Sink_1}$ be the sets of facilities selected using the above listed four approaches in order,  $t_{MILP}$, $t_{LP}$, $t_{Sink_2}$ and $t_{Sink_1}$ be the solution times, and $J_{MILP}$, $J_{LP}$, $J_{Sink_2}$ and $J_{Sink_1}$ be the respective objective values. For different values of $k$, we assess the quality of our algorithm in terms of the following metrics: (i) the percentage overlaps on the final selected set of facilities, (ii) the time taken to select the all the facilities, (iii) the relative differences on the original objective function in~\eqref{eqn:set_selection1:objective}, and (iv) the approximation of the demand allocation problem $\mathbb{P}_3$ by the multi-stage Sinkhorn oracle relative to  standard LP solver.

In Fig. \ref{fig:compare}(b) we plot the percentage of facilities selected using the three oracles that are common with the solution from MILP, computed as: $\calS_{algo}^{MILP} = \frac{|\calS_{algo} \cap \calS_{MILP}|}{|\calS_{MILP}|} \times 100 \%$, $\forall algo \in \{LP, Sink_2, Sink_1\}$. For most values of $k$ the overlap is greater than 80$\%$. In Fig. \ref{fig:compare}(c) we plot the relative difference in objective value against the MILP solution calculated as: $\Delta J_{algo}^{MILP} = \frac{J_{algo}-J_{MILP}}{J_{MILP}}\times 100\%$, $\forall algo \in \{LP, Sink_2, Sink_1\}$), and find that both the LP and the multi-stage Sinkhorn oracles are within \textbf{3}$\%$. Any (acceptable) loss in optimality from the Sinkhorn oracles is compensated by the significant gain in computational time to select facilities as seen in Fig.~\ref{fig:compare}(a). The execution time of the multi-stage Sinkhorn oracle can be further reduced by running the stage-2 OT instances in parallel, and also by making parallel calls to the value oracle for multiple facility options. Such a study is beyond the scope of this paper. As expected, the \textit{single-stage} Sinkhorn oracle is faster than its multi-stage counterpart but also has a higher variation from the optimal objective value, as shown in Fig.~\ref{fig:compare}(c).

We also make a stand-alone comparison of the proposed multi-stage Sinkhorn algorithm with a standard commercial LP solver, without any incremental facility selection algorithm. We execute the problem $\mathbb{P}_3$ with the selected set of facilities $\calS_{Sink_2}$, for different values of cardinality constraint $k$, using: (i) an LP solver and, (ii) the multi-stage Sinkhorn algorithm~\ref{algo:dual_sinkhorn}. Let the solution objectives be represented as $g_{LP}(\calS_{Sink_2})$ and $g_{Sink_2}(\calS_{Sink_2})$ respectively. The relative difference, $\Delta g^{LP}_{Sink_2} = \frac{g_{Sink_2}(\calS_{Sink_2})-g_{LP}(\calS_{Sink_2})}{g_{LP}(\calS_{Sink_2})}\times 100\%$, between the objective from multi-stage Sinkhorn and the optimal value from LP solver is less than $1\%$, as shown in Fig.~\ref{fig:compare}(d). We thus obtain at least \textbf{20} fold speed up \commentbyKG{just in execution of \emph{one instance of $\mathbb{P}_3$} }by replacing the LP solver with multi-stage Sinkhorn as the value oracle, without any significant compromise in  the quality of the two results.

\begin{figure}[t]
  \centering
  \includegraphics[width=0.9\linewidth,trim=0.2cm 0cm 0cm 0cm]{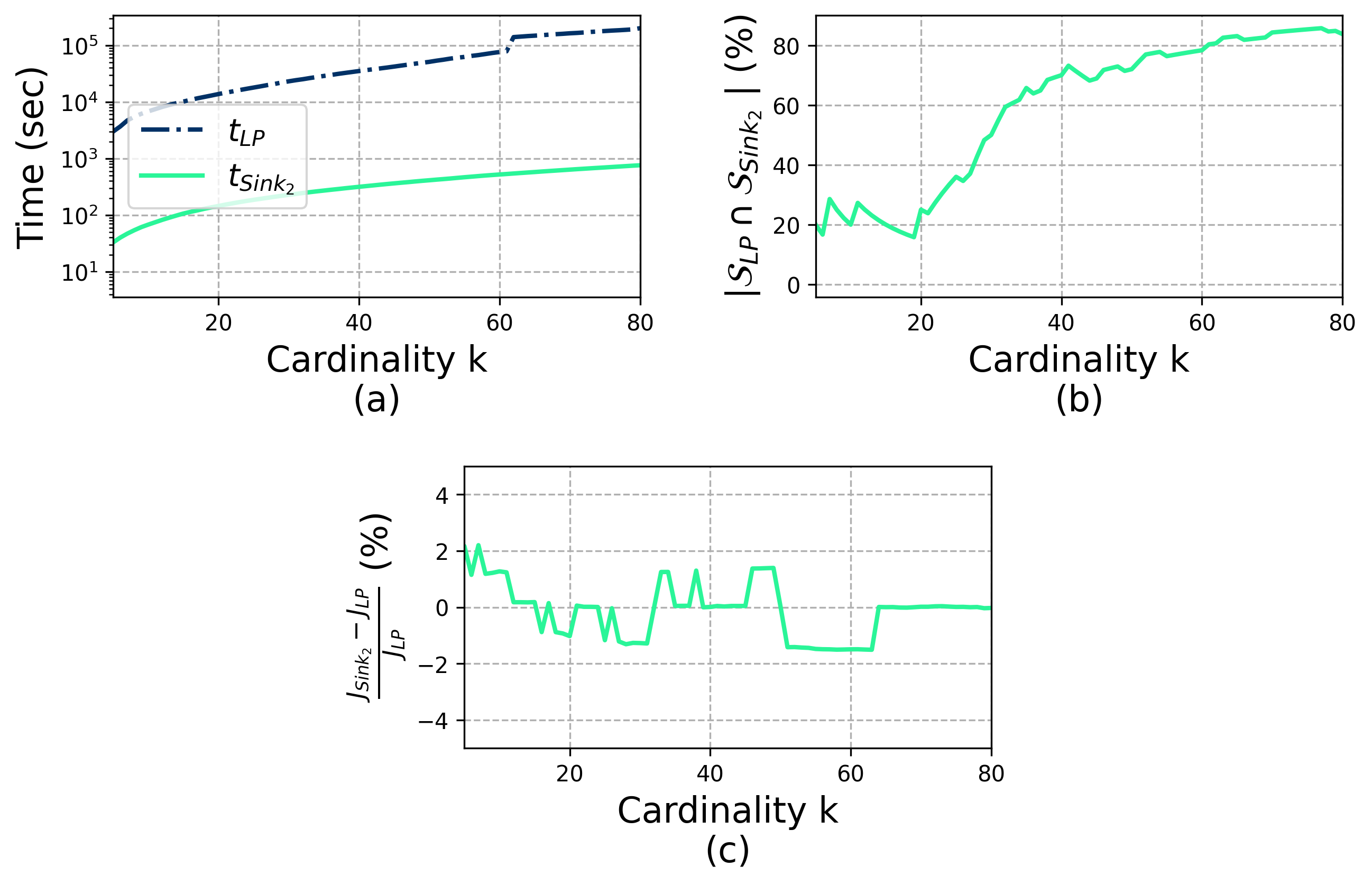}
  \caption{Results for the large network with $n=10000$, $m=400$, $|\calE|>10^6$ . (a) Time taken to select all the facilities, (b) Percentage overlap between the final set of selected facilities, and (c) Relative difference between the objective values.}
  \label{fig:compare_large_netw}
\end{figure}

Lastly, we stress test the scalability of our approach by selecting facilities for a much larger network with 400 facilities, 10000 clients, and \emph{more than 10 million flow paths} with three different shipment channels. It takes a tremendously large amount of time to solve an MILP for this scale of the FLP and is practically infeasible to be used even for comparison studies. For e.g., the best possible solution obtained using MILP in 48 hrs for $k=80$ has an objective value which is worse by 40$\%$ than our approximate solution. We therefore use the stochastic distorted greedy algorithm with an LP value oracle to benchmark our solution in this case. The total time taken to select the facilities for different cardinality values $k$ is shown in Fig.~\ref{fig:compare_large_netw}(a). While the direct LP solver took more than $200000$ secs to identify $80$ facilities, the multi-stage Sinkhorn oracle was able to choose the same in $\approx750$ secs. The overlap between the solutions from using LP and multi-stage Sinkhorn as value oracles is shown in Fig. \ref{fig:compare_large_netw}(b). While this overlap is less for lower values of $k$, the difference in objective values from the two solutions never exceeds $2.5\%$ as seen in Fig.~\ref{fig:compare_large_netw}(c). In fact, since the subsets $\mathcal{M}^r \subseteq \mathcal{M}-\mathcal{S}$ in each iteration of the stochastic distorted greedy algorithm are chosen randomly, for some $k$ values the multi-stage Sinkhorn value oracle even does better than the LP solver. This re-establishes the significance of our method.

\section{Conclusion}
We introduced the novel multi-channel capacitated facility location problem, where the facilities can service the demand at the client locations using different fulfilment channels, each with their own service cost and channel capacities. As opposed to a standard FLP, this creates complex coupling constraints between the shipment channels. Though polynomial time iterative algorithms exists for the capacitated FLP, we argued that these methods do not scale as: (i) they perform complicated operations allowing to open or close a subset of facilities relative to the current solution which are computationally intensive, (ii) rely on using the LP solvers as the value oracle for determining the cost of demand allocation given the set of facilities $\calS$, to which repeated invocations are made over multiple iterations. We proposed a scalable alternative by: (i) modelling the FLP as a submodular maximisation problem and using the well-known greedy technique which only performs the relatively simpler operation of incrementally adding a new facility to the current selection, (ii) replacing the direct LP solver in the value oracle by a multi-stage instances OT problem, each solvable in lightspeed using Sinkhorn iterations. The theoretical and convergence guarantees of these greedy and Sinkhorn approaches directly apply to our framework. As observed in our experiments, our OT based algorithm leads to a $\mathbf{20}$ fold decrease in computational time over state-of-the-art commercial LP solvers and over $\mathbf{100}$ fold decrease in the same as compared to the direct MILP solvers, with only a marginal variation of less than $\mathbf{3\%}$ from the optimal objective value. 

To the best of our knowledge, all the variants of FLP studied thus far are \emph{static} in nature, where a one-time decision is made to open a facility given the various costs, capacity, and the demand values. However in the real world, client demands, facility capacities, service cost, etc. constantly vary over time. A straightforward approach for this \emph{dynamic} setting would be to be to treat each time instance as a standalone static problem, and freshly determine the set of facilities to open for that time period. A number of operational constraints would preclude continuously opening and shutting down facilities. For instance, a supply chain system would not prefer to open a facility in the first year, close it down in the following year only to open it again. As our future work, we would like to investigate the multi-year, dynamic variant of FLP with changing input parameters over the years, where one could have the option to open a facility any time during the multi-year period, but constrained from closing it subsequently.

\bibliographystyle{splncs04}
\bibliography{facilityrecommender}
\newpage
\appendix
\section{Facility Location Problem (MILP)}
\label{sec:MILP}
Given the supply chain network described in Section~\ref{sec:problem}, finding $\mathcal{S}^* = \{i  | y_i = 1, i \in \mathcal{M}\}$ where the value of $y_i$,  $\forall i \in \mathcal{M}$ is obtained by solving the below problem,
\begin{equation}\label{eqn:MILP:objective}
 \mathbb{P}_6: \min_{x_{i,j,e}, y_{i}}  J  = \sum_{i \in \mathcal{M}} F_{i}  y_{i} + \sum_{(i,j,e) \in \mathcal{E}} c_{i,j,e}  x_{i,j,e} d_{j} + K \left(\sum_{j \in \mathcal{N}} d_{j} - \sum_{(i,j,e) \in \mathcal{E}} x_{i,j,e} d_{j}\right) 
\end{equation}
such that,
\begin{align}
\sum_{i \in \mathcal{M}} \sum_{e \in E} x_{i,j,e} & \leq 1 & \quad \forall j \in \mathcal{N} \label{eqn:MILP:demand} \\
\sum_{j \in \mathcal{N}} \sum_{e \in E} x_{i,j,e} d_j & \leq fcap_i & \quad \forall i \in \mathcal{M} \label{eqn:MILP:fcap}\\
\sum_{j \in \mathcal{N}} x_{i,j,e} d_j &\leq ccap_{i,e} & \quad \forall i \in \mathcal{M}, \quad \forall e \in E \label{eqn:MILP:ccap}\\
 \sum_{i \in \mathcal{M}} y_i & \leq k & \label{eqn:MILP:cardinality}\\
  x_{i,j,e} & \leq y_i & \quad \forall i \in \mathcal{M}, \quad \forall e \in E, \quad \forall j \in \mathcal{N} \label{eqn:MILP:no_ship_from_closed_facility}\\
 x_{i,j,e} &\geq 0 & \quad \forall i \in \mathcal{M}, \quad \forall e \in E, \quad \forall j \in \mathcal{N} \label{eqn:MILP:positive_ship}\\
  x_{i,j,e} &\leq 0 & \quad \forall (i, j, e) \in \{\mathcal{M} \times \mathcal{N} \times E\}\backslash \mathcal{E} \label{eqn:MILP:negative_ship}\\
   y_i & \in \{0,1\} & \quad \forall i \in \mathcal{M}. \label{eqn:MILP:binary_variable}
  \end{align}
  
The above problem is an MILP model of the FLP under consideration. Each facility $i$ is modelled as a Boolean variable $y_i$, taking values in the set $\{0,1\}$. The value $y_i=1$ indicates that the corresponding facility is opened and closed otherwise. The value $x_{i,j,e}$ is the proportion of demand $d_j$ at the $j^{th}$ client fulfilled by the $i^{th}$ facility using the shipment channel $e$. The constraint in equation~\eqref{eqn:MILP:no_ship_from_closed_facility} guarantees that orders are fulfilled using the opened facilities only, while equations~\eqref{eqn:MILP:positive_ship} and \eqref{eqn:MILP:negative_ship} ensure that the network connectivity constraints are followed. The cardinality constraint on the total number of selected facilities is imposed in \eqref{eqn:MILP:cardinality}. The optimization objective given in~\eqref{eqn:MILP:objective} of the FLP has three components: (i) first term is the fixed cost to open the selected facilities, (ii) the middle term denotes the total cost to ship, and (iii) lastly a penalty on any unfulfilled demand. The parameter $K > \max\limits_{(i,j,e) \in \mathcal{E}} c_{i,j,e}$ controls the trade-off between the fulfilled and unfulfilled demand. Sometimes an FLP solution may not be able to meet all the demand either due to insufficient capacity, or insufficient supply due to transportation restrictions (network constraints). Allowing a slack on the fulfilled demand in~\eqref{eqn:MILP:demand} ensures the feasibility of the problem.

\section{Proofs}
\label{sec:prooofs}
\subsection{Proof of Lemmma~\ref{lemma:decouple1}}
Consider $\sum\limits_{e \in E_i} ccap_{i,e} \leq fcap_i$. Summing up equation~\eqref{eqn:set_selection1:ccap} over the channels $e \in E_i$, we get
\begin{equation*}
    \sum\limits_{e \in E_i} \sum_{j \in \mathcal{N}_i}  x_{i,j,e} d_j \leq \sum\limits_{e \in E_i} ccap_{i,e} \leq fcap_i.
\end{equation*}
and the condition~\ref{eqn:set_selection1:fcap} is trivially satisfied.

\subsection{Proof of Lemma~\ref{lemma:decouple2}}
Let $a_{i,e} \leq ccap_{i,e}$ represent the utilisation of the channel capacity $e$ from the facility $i$ in the solution to problem $\mathbb{P}_3$. It can be expressed as: $a_{i,e} = \sum\limits_{j \in \mathcal{N}_i} x^{\ast}_{i,j,e} d_j$, where $x^{\ast}_{i,j,e}$ is the global optimiser for $\mathbb{P}_3$.
If $a_{i,e^{\prime}} > fcap_i - \left(\sum\limits_{e \in E_i-\{e^{\prime}\}} ccap_{i,e}\right) \geq 0$, then exits a channel $e \not= e^{\prime}$ such that 
$a_{i,e} < ccap_{i,e}$, otherwise the facility capacity constraint in  will be violated. Consider any client $j \in \mathcal{N}_i$ that is serviced by $e^{\prime}$. 
As $p_{i,j,e^{\prime}} \leq p_{i,j,e}$ by the condition $e^{\prime} = \argmin\limits_{e \in E_i} p_{i,j,e},$ $\forall j \in \mathcal{N}_i$ of Lemma~\ref{lemma:decouple2}, we could move a $\delta > 0$ allocation from $e^{\prime}$ to $e$ for the client $j$, i.e. $\hat{x}_{i,j,e^{\prime}}  = x^{\ast}_{i,j,e^{\prime}} -\delta$, and $ \hat{x}_{i,j,e} = x^{\ast}_{i,j,e} + \delta$ and overall objective value will either increase (resulting in a contradiction) or stay the same. In either case, $a_{i,e^{\prime}}$ should or can be made to satisfy the condition $a_{i,e^{\prime}} \leq fcap_i - \left(\sum_{e \in E_i-\{e^{\prime}\}} ccap_{i,e}\right)$, and the capacity of $ccap_{i,e^{\prime}} - a_{i,e^{\prime}}$ will remain not utilised. The capacity of channel $e^{\prime}$ can hence be reduced without affecting the global optimum.

\section{Channel Decoupling Example}
\label{sec:channeldecouplingexample}

\begin{figure}[!t] 
  \centering
  \includegraphics[width=0.5\linewidth]{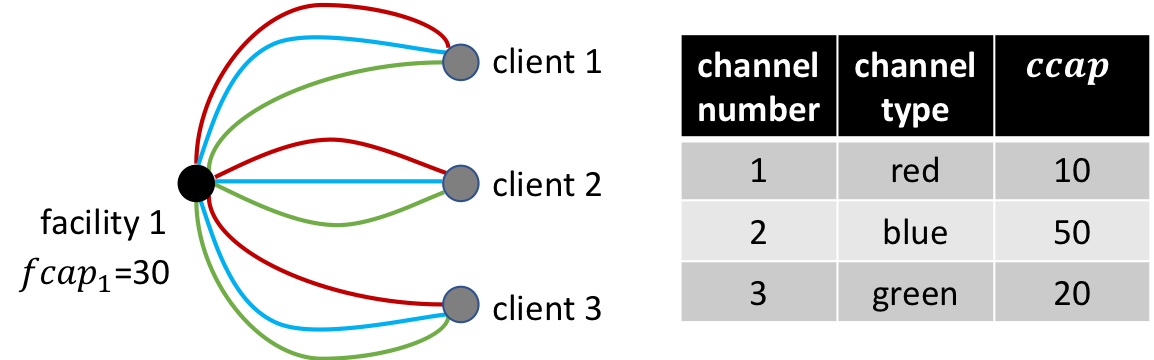}
  \caption{An example for channel decoupling}
  \label{fig:channel_decoupling}
\end{figure}

Using an example, we explain Lemmas~\ref{lemma:decouple1} and \ref{lemma:decouple2} that present an approach to decouple the outflow on multiple shipment channels under favourable conditions. Consider a sample single facility sub-network $\mathbf{T}_{1}$ in Fig.~\ref{fig:channel_decoupling}. Facility-1 is connected to three clients through three different fulfilment channels indicated in red, blue, and green. The facility capacity and the channel capacities satisfy $fcap_1 < ccap_{1,1} + ccap_{1,2} + ccap_{1,3}$. The profit numbers are given in Table~\ref{table:channel_decoupling}. We discuss multiple scenarios:\\\\
\textbf{Case 1:} Based on the case-1 profit numbers in Table~\ref{table:channel_decoupling}, red (channel-1) is ubiquitously the least profitable channel at this facility, whose channel capacity can be reduced. The new channel capacity becomes $
ccap_{1,1} = \max \left(0, fcap_1 - \left( ccap_{1,2} + ccap_{1,3}\right) \right) = 0$,
and this zero capacity channel can be completely removed from facility-1. Even with this capacity update, $fcap_1 < ccap_{1,2} + ccap_{1,3}$. We can re-apply Lemma~\ref{lemma:decouple2} on the next least profitable channel, blue, across all clients and reduce its channel capacity to $ccap_{1,2} =$ $ \max (0, fcap_1 - $ $ccap_{1,3}) = 10$.
These updated channel capacities satisfy Lemma~\ref{lemma:decouple1}, and outflow on the channels is decoupled.\\\\
\textbf{Case 2:} As seen in Table~\ref{table:channel_decoupling}, red is still the least profitable channel. Its channel capacity can be reduced to zero after application of Lemma~\ref{lemma:decouple2}. But, unlike case-1, there is no channel which is the second-least profitable from facility-1 to \emph{all} its clients. Therefore, a second application of Lemma~\ref{lemma:decouple2} is not possible, and the blue and green channels remain coupled through $fcap_1$.\\\\
\textbf{Case 3:} This case represents a very important scenario where
 red is the least profitable \textit{available} channel at facility-1, but since blue channel (channel-2) is unavailable to service client-3, $p_{1,3,2} = 0$. The condition of Lemma~\ref{lemma:decouple2} is not satisfied, and shipment channels cannot be decoupled.
\begin{table}
    \centering
    \caption{Profit to ship unit item for network in Fig. \ref{fig:channel_decoupling}}
    \label{table:channel_decoupling}
    \begin{adjustbox}{width=0.5\columnwidth,center}
		\begin{tabular}{l|llll}
			\hline
			 Case & channel type & client 1 & client 2 & client 3 \\
			\hline
			      & red   & 3.5  & 4.2  &   3.9   \\
			Case 1  & blue  & 4.5  & 4.3   &   4.1   \\
                    & green & 4.9  & 4.8   &   4.7   \\
			\hline
				    & red   & 3.5  & 0   &   3.9   \\
			Case 2  & blue  & 3.5  & 4.4   &   4.1   \\
                    & green & 4.9  & 3.6   &   4.7  \\
			\hline
				    & red   & 3.5  &  4.2  &   3.9   \\
			Case 3  & blue  & 4.5  & 4.3   &   0   \\
                    & green & 4.9  & 4.8   &   4.7   \\
			\hline
		\end{tabular}
  \end{adjustbox}
\end{table}
\end{document}